\documentclass{amsart}

\usepackage{amsmath,amssymb}

\theoremstyle{plain}

\newtheorem{theorem}{Theorem}[section]
\newtheorem{proposition}[theorem]{Proposition}
\newtheorem{lemma}[theorem]{Lemma}
\newtheorem{corollary}[theorem]{Corollary}
\newtheorem{conjecture}[theorem]{Conjecture}

\theoremstyle{definition}

\newtheorem{remark}[theorem]{Remark}

\newcommand{\C}{\mathcal C}
\newcommand{\G}{\mathcal G}

\newcommand{\p}{\mathfrak p}
\renewcommand{\O}{\mathcal O}
\newcommand{\Q}{\mathbf Q}
\renewcommand{\P}{\mathcal P}
\newcommand{\Z}{\mathbf Z}

\newcommand{\dmn}{\delta_{m}^{n}}

\newcommand{\pib}{\bar{\pi}}
\newcommand{\pn}{\varphi(n)}

\newcommand{\inj}{\hookrightarrow}

\newcommand{\Ct}[2]{\C_{\sqrt{#1}}^{#2}}
\newcommand{\Ctt}[3]{\C_{\sqrt[#1]{#2}}^{#3}}
\newcommand{\ls}[2]{\left( \frac{#1}{#2} \right)}
\renewcommand{\mod}[1]{\,(\text{mod }#1)}

\renewcommand{\th}{\text{th}}

\DeclareMathOperator{\Gal}{Gal}

\begin{document}

\title[Power residues]{Power residues of Fourier coefficients of elliptic
curves with complex multiplication}
\author{Tom Weston and Elena Zaurova}
\address[Tom Weston]{Dept.\ of Mathematics, University of
Massachusetts, Amherst, MA}
\address[Elena Zaurova]{Dept.\ of Mathematics, University of Massachusetts,
Amherst, MA}

\email[Tom Weston]{weston@math.umass.edu}
\email[Elena Zaurova]{ezaurova@student.umass.edu}

\thanks{The first author was supported by NSF grant DMS-0440708; the second
author was supported by generous contributions from alumni of the University
of Massachusetts}

\maketitle

Let $E$ be an elliptic curve over $\Q$.  For any $m \geq 1$
and set of primes $\C$ (contained in the set of primes congruent to one
modulo $m$) we define $\delta_{m}^{1}(E;\C)$ as
the relative density (in the set of $p \in \C$ which are ordinary for $E$)
of primes $p \in \C$ for which the $p^{\th}$ Fourier coefficient of
$E$ is an $m^{\th}$-power modulo $p$.
In \cite{prw} it was
conjectured that $\delta_{m}^{1}(E;\C) = \frac{1}{m}$ whenever $E$
does not have complex multiplication and $\C$ is a set of primes defined
by Galois theoretic conditions.  In the present paper we extend these
conjectures to the case of elliptic curves with complex multiplication; we
also prove our conjectures for certain small values of $m$.

To be more precise, fix an imaginary quadratic field $K$ of class number one
and let $E$ denote an elliptic curve with complex multiplication by the
ring of integers $\O_{K}$; we write $w$ for the order of $\O_{K}^{\times}$.
For any divisor $n$ of $m$ we consider the density $\delta_{m}^{n}(E;\C)$
of $p \in \C$ for which the $m^{\th}$ power residue symbol of the 
$p^{\th}$ Fourier coefficient of $E$ modulo $p$ is a primitive 
$n^{\th}$ root of unity.
We compute the density $\delta_{m}^{n}(E;\C)$ (in terms of certain simpler
densities) for any $m$ dividing $w$; most of these computations were
essentially done in \cite{prw}, with the exception of $K=\Q(i)$ and $m=4$
(which is significantly more involved).
These densities are often different from the naive expectation
$\frac{\varphi(n)}{m}$.  For general $m$, we conjecture that the
density $\delta_{m}^{n}(E;\C)$ differs from
$\frac{\varphi(n)}{m}$ only to the extent that such a difference is forced
upon it by its relation to the known density $\delta_{m'}^{n'}(E;\C)$ with
$m'=(m,w)$ and $n'$ an appropriate divisor of $m'$.  We make
our conjecture entirely explicit in the case that $\C$ consists of all
primes congruent to one modulo $m$.

We now outline the contents of the paper.  In Section 1 we set our notation
for densities and give the basic density computation coming from the Chebotarev
density theorem.  We give careful statements of the known results for $m$
dividing $w$ in Section 2; this is done most efficiently by regarding
all elliptic curves with complex multiplication by $\O_{K}$ as twists
of a fixed such curve.  We also give some preliminary density
computations for later
use in explicating our general conjectures.  Those conjectures are stated
in Section 3, where we also verify certain natural compatibilities.  In Section
4 we give the proof, based on biquadratic reciprocity, 
of our conjecture in the case $K=\Q(i)$ and $m=4$

We emphasize that despite the essentially elementary nature of our approach
for $m$ dividing $w$, it remains our opinion that entirely new methods will
be required to approach the general case.

It is a pleasure to thank Farshid Hajir for numerous helpful conversations.

\section{Densities}

\subsection{Preliminaries}

For a prime ideal $\p$ of a finite extension $K$ of $\Q(\zeta_{m})$ and
$\alpha \in \O_{K} - \p$, we write
$\ls{\alpha}{\p}_{m}$ for the $m^{\th}$ power residue symbol of $\alpha$
modulo $\p$.  Thus $\ls{\alpha}{\p}_{m} \in \mu_{m}$ and
$$\ls{\alpha}{\p}_{m} \equiv \alpha^{\frac{N(\p)-1}{m}} \pmod{\p}.$$
Most often we will be interested only in the order of
$\ls{\alpha}{\p}_{m}$ and in the case that $\alpha \in \Z$ and
$K=\Q(\zeta_{m})$.  When this is the case, by abuse of notation we
simply write $\ls{\alpha}{p}_{m}$ to mean $\ls{\alpha}{\p}_{m}$ for some
prime ideal $\p$ of $\Q(\zeta_{m})$ above $p$.  We emphasize that while
the precise root of unity $\ls{\alpha}{p}_{m}$ is not well-defined, its
order is well-defined.

For a set of positive rational
primes $\P$ we define the {\it zeta function} $\zeta(s;\P)$ of $\P$ by
$$\zeta(s;\P) = \sum_{p \in \P} p^{-s};$$
this converges for $s > 1$.  If $\P$ and $\P'$ are sets of
primes with $\P$ of positive density (in the sense that
$\lim_{s \to 1^{+}} \zeta(s;\P)$ diverges), 
we define the {\it  relative density} 
$\rho_{\P}(\P')$ of $\P'$ in $\P$ by
$$\rho_{\P}(\P') = \underset{s \to 1^{+}}{\lim} \frac{\zeta(s;\P \cap \P')}
{\zeta(s;\P)}$$
assuming it exists.  We note that
$$\rho_{\P}(\P') = \underset{x \to \infty}{\lim} \frac{\# \{p \in \P \cap \P' \,;\, 
p < x \}}{\# \{p \in \P \,;\, p < x \}}$$
when the latter limit exists (which it will in all cases we consider).

We now introduce the sets of primes we will work with.  For a finite Galois
extension $K/\Q$ and a union $C$ of conjugacy classes in
$\Gal(K/\Q)$, we write $\C_{K}^{C}$ for the set of rational primes $p$,
unramified in $K/\Q$, with Frobenius over $K$ lying in $C$.  Recall that by
the Chebotarev density theorem the set of primes $\C_{K}^{C}$ has absolute
density equal to $\frac{\#C}{[K:\Q]}$.
We say that
a set $\P$ of primes is {\it Chebotarev} if it agrees with some
$\C_{K}^{C}$ up to finite sets.

The basic density result we will need is the following.

\begin{lemma} \label{lemma:dens}
Let $K_{1},K_{2}$ be finite Galois extensions of $\Q$ and fix subsets
$C_{i} \subseteq \Gal(K_{i}/\Q)$ stable under conjugation.  Then
$$\rho_{\C_{K_{1}}^{C_{1}}}(\C_{K_{2}}^{C_{2}}) = 
\frac{\# \{ (\sigma_{1},\sigma_{2}) \subseteq
C_{1} \times C_{2} \,;\, \sigma_{1}|_{K_{1} \cap K_{2}} =
\sigma_{2}|_{K_{1} \cap K_{2}} \}}{\# C_{1} \cdot [K_{2}:K_{1} \cap K_{2}]}.$$
\end{lemma}
\begin{proof}
Let $\pi_{i} : \Gal(K_{1}K_{2}/\Q) \to \Gal(K_{i}/\Q)$ denote the
natural surjection.  We then have
$$\C_{K_{i}}^{C_{i}} = \C_{K_{1}K_{2}}^{\pi_{i}^{-1}(C_{i})}.$$
In particular,
$$\C_{K_{1}}^{C_{1}} \cap \C_{K_{2}}^{C_{2}} = \C_{K_{1}K_{2}}^{
\pi_{1}^{-1}(C_{1}) \cap \pi_{2}^{-1}(C_{2})}.$$
By the Chebotarev density theorem we thus have
$$\rho_{\C_{K_{1}}^{C_{1}}}(\C_{K_{2}}^{C_{2}}) =
\frac{\# \pi_{1}^{-1}(C_{1}) \cap \pi_{2}^{-1}(C_{2})}
{\# \pi_{1}^{-1}(C_{1})}.$$
As 
$$\# \pi_{1}^{-1}(C_{1}) = [K_{1}K_{2}:K_{1}] \cdot \# C_{1} =
[K_{2}:K_{1} \cap K_{2}] \cdot \# C_{1},$$
the lemma thus follows from noting that the injection
$$\pi_{1} \times \pi_{2} : 
\pi_{1}^{-1}(C_{1}) \cap \pi_{2}^{-1}(C_{2}) \inj C_{1} \times C_{2}$$
has image
$$\{ (\sigma_{1},\sigma_{2}) \subseteq
C_{1} \times C_{2} \,;\, \sigma_{1}|_{K_{1} \cap K_{2}} =
\sigma_{2}|_{K_{1} \cap K_{2}} \}.$$
\end{proof}

We fix notation for the most common Chebotarev sets we will encounter.  For
relatively prime integers $a$ and $b$ we write $\C_{b}^{a}$ for the set
of primes congruent to $a$ modulo $b$; it is a Chebotarev set with
$K=\Q(\mu_{b})$.  For $t \in \Q^{\times}$, $m \geq 1$ and $\zeta \in \mu_{m}$
we write $\Ctt{m}{t}{\zeta}$ for the Chebotarev set of primes
$p \equiv 1 \mod{m}$ such that $\ls{t}{p}_{m}$ is conjugate to $\zeta$.
We often simply write $\Ctt{m}{t}{+}$ (resp.\ $\Ctt{m}{t}{-}$) for
$\Ctt{m}{t}{1}$ (resp.\ $\Ctt{m}{t}{-1}$).  For example,
$\Ct{t}{+}$ (resp.\ $\Ct{t}{-}$) is the set of odd primes which are split
(resp.\ inert) in $\Q(\sqrt{t})$, where we consider all primes to
be split in $\Q$ in the case that $t \in \Q^{\times 2}$.

\subsection{Elliptic curves}

Fix an elliptic curve $E$ over $\Q$ and $m \geq 1$.  For any set of primes
$\C$ we define
$$\P_{m}(E;\C) = \bigl\{ p \in \C_{m}^{1} \cap \C \,;\, a_{p}(E) \not\equiv 0
\mod{p} \bigr\}.$$
When $\P_{m}(E;\C)$ has positive density,
for $n$ dividing $m$ we define $\dmn(E;\C)$ as the relative density of
$$\P_{m}^{n}(E;\C) = \left\{ p \in \P_{m}(E;\C) \,;\,
\ls{a_{p}(E)}{p}_{m} \text{~has exact order $n$} \right\}$$
in $\P_{m}(E;\C)$.  

In \cite{prw} it was conjectured that if $E$ does not have complex
multiplication, then $\delta_{m}^{1}(E;\C) = \frac{1}{m}$ for any 
Chebotarev set $\C$ contained in $\C_{m}^{1}$; more generally, one
expects that $\dmn(E;\C) = \frac{\pn}{m}$ for any $n$ dividing $m$.
Our goal in this paper is to formulate (and prove for small $m$)
analogous conjectures when $E$ does have complex multiplication.

\section{Elliptic curves with complex multiplication}

\subsection{Twisting}

Fix a discriminant
$$d \in \{-3,-4,-7,-8,-11,-19,-43,-67,-163 \}.$$  Let $w_{d}$ denote the
number of units in the ring of integers of $\Q(\sqrt{d})$.
For $d \neq -3,-4,-8$, let $E^{d}$ denote an elliptic curve of conductor
$d^{2}$ with complex multiplication by $\Z[\frac{1+\sqrt{d}}{2}]$;
the curve $E^{d}$ is determined up to isogeny, which will suffice for our
purposes.  We let $E^{-8}$ denote an elliptic curve in the
isogeny class 256D in \cite{Cremona}, with complex multiplication by
$\Z[\sqrt{-2}]$.

For $t \in \Q^{\times}$ we define an elliptic curve $E^{d}_{t}$ as follows:
\begin{itemize}
\item Let $E^{-3}_{t}$ denote the elliptic curve with
Weierstrass equation $y^{2} = x^{3}+16t$.
\item Let $E^{-4}_{t}$ denote the elliptic curve with Weierstrass
equation $y^{2}=x^{3}-tx$.
\item Let $E^{-7}_{t}$ denote the quadratic twist of $E^{-7}$ by $-t$.
\item For $d \leq -8$, let $E^{d}_{t}$ denote the quadratic twist of $E^{d}$ by
$t$.
\end{itemize}
(The slightly different twist in the case $d=-7$ is necessary because
$2$ is split in $\Q(\sqrt{-7})$.)
In particular, for $d \neq -3,-4$ we have
\begin{equation} \label{eq:tw}
a_{p}(E^{d}_{t}) = \begin{cases} \ls{t}{p}a_{p}(E^{d}_{1}) &
p \in \Ct{d}{+} \text{~and~} (t,p)=1; \\
0 & p \in \Ct{d}{-}. \end{cases}
\end{equation}
In any case, for any Chebotarev set $\C$ the set of primes
$\P_{m}(E^{d}_{t};\C)$ differs from
$\C \cap \C_{m}^{1} \cap \Ct{d}{+}$ by a finite set.  We thus can,
and for the remainder of this section will, assume that
$\C \subseteq \C_{m}^{1} \cap \Ct{d}{+}$.

\subsection{Densities for $m \mid w_{d}$}

Let $w_{d}$ 
denote the order of the group of units $\O_{\Q(\sqrt{d})}^{\times}$;
thus $w_{-3} = 6$, $w_{-4} = 4$ and $w_{d} = 2$ for $d \neq -3,-4$.
We now recall the known formulae for $\delta_{m}^{n}(E^{d}_{t};\C)$
for all $m$ dividing $w_{d}$.  We begin with the case $m=2$ when such results
exist for all $d$.

\begin{proposition} \label{prop:m2}
Let $\C$ be a Chebotarev set contained in $\Ct{d}{+}$.
Then
$$\delta_{2}^{1}(E^{d}_{t};\C) = \begin{cases}
\rho_{\C}(\Ctt{4}{d}{+}) + \rho_{\C}(\C_{4}^{3} \cap \Ct{t}{-}) &
d \neq -4; \\
\rho_{\C}(\C_{8}^{1}) + \rho_{\C}(\C_{8}^{5} \cap \Ct{t}{-}) & d = -4.
\end{cases}$$
$$\delta_{2}^{2}(E^{d}_{t};\C) = \begin{cases}
\rho_{\C}(\Ctt{4}{d}{-}) + \rho_{\C}(\C_{4}^{3} \cap \Ct{t}{+}) &
d \neq -4; \\
\rho_{\C}(\C_{8}^{5} \cap \Ct{t}{+}) & d = -4;
\end{cases}$$
for any $t \in \Q^{\times}$.
\end{proposition}
\begin{proof}
This is immediate from the formulae
$$\ls{a_{p}(E^{d}_{t})}{p} = \begin{cases}
\ls{-d}{p}_{4} & p \in \C_{4}^{1} \cap \Ct{d}{+}; \\
-\ls{t}{p} & p \in \C_{4}^{3} \cap \Ct{d}{+};
\end{cases}$$
(for $d \neq -4$) and
$$\ls{a_{p}(E^{-4}_{t})}{p} = \begin{cases}
1 & p \in \C_{8}^{1}; \\
-\ls{t}{p} & p \in \C_{8}^{5};
\end{cases}$$
of \cite{prw}.  (Note that $\ls{-d}{p}_{4} = \pm 1$ for
$p \in \Ct{d}{+}$.)
\end{proof}

Before we can state our result for $d=-4$ and $m=4$ we must introduce
some notation.
For $\alpha \in \{1,1+4i,5,5+4i\}$, let $\G_{8}^{\alpha}$ denote
the set of rational primes $p \equiv 1 \mod{8}$ for which one
(or equivalently both) of the primary divisors of $p$ in
$\Q(i)$ are congruent to $\alpha$ modulo $8$.  (Recall that an element
$\alpha \in \Z[i]$ is said to be {\it primary} if $\alpha \equiv 1
\mod{2+2i}$.)  These are Chebotarev sets for the ray class field
$K=\Q(\zeta_{16},\sqrt[4]{2})$ of $\Q(i)$ of conductor $8$.  These four
sets partition $\C_{8}^{1}$.

The next proposition is proved in Section~\ref{sec:gi}.

\begin{proposition} \label{prop:m4}
Let $\C$ be a Chebotarev set contained in
$\C_{4}^{1}$.  Then
\begin{align*}
\delta_{4}^{1}(E^{-4}_{t};\C) &=
\rho_{\C}(\G_{8}^{1}) + \rho_{\C}(\G_{8}^{5} \cap \Ct{t}{-}) +
\rho_{\C}(\G_{8}^{5+4i} \cap \Ct{t}{+}) + 
\rho_{\C}(\C_{8}^{5} \cap \Ctt{4}{t/2}{+}); \\
\delta_{4}^{2}(E^{-4}_{t};\C) &=
\rho_{\C}(\G_{8}^{1+4i}) + \rho_{\C}(\G_{8}^{5} \cap \Ct{t}{+}) +
\rho_{\C}(\G_{8}^{5+4i} \cap \Ct{t}{-}) + 
\rho_{\C}(\C_{8}^{5} \cap \Ctt{4}{t/2}{-}); \\
\delta_{4}^{4}(E^{-4}_{t};\C) &=
\rho_{\C}(\C_{8}^{5} \cap \Ct{t}{+});
\end{align*}
for any $t \in \Q^{\times}$.
\end{proposition}

Finally, when $d=-3$ and $m=3$ or $m=6$ we have the following result.
Let $\omega$ denote a primitive third root of unity.

\begin{proposition} \label{prop:m36}
Let $\C$ be a Chebotarev set contained in
$\C_{3}^{1}$.  Then
\begin{align*}
\delta_{3}^{1}(E^{-3}_{t};\C) &=
\rho_{\C}(\C_{9}^{1}) + \rho_{\C}(\C_{9}^{4,7} \cap \Ctt{3}{t}{+}); \\
\delta_{3}^{1}(E^{-3}_{t};\C) &=
\rho_{\C}(\C_{9}^{4,7} \cap \Ctt{3}{t}{\omega});
\end{align*}
for any $t \in \Q^{\times}$.  Furthermore,
$$\delta_{6}^{n}(E^{-3}_{t};\C) = \delta_{2}^{(n,2)}(E^{-3}_{t};\C) \cdot
\delta_{3}^{(n,3)}(E^{-3}_{t};\C)$$
for any $t \in \Q^{\times}$.
\end{proposition}

Here by $\C_{9}^{4,7}$ we of course mean the set of all primes $p$
which are equivalent to either $4$ or $7$ modulo $9$.

\begin{proof}
The densities for $m=3$ follow from the formula
$$\ls{a_{p}(E^{-3}_{t})}{\pi}_{3} = \begin{cases} 1 & p \equiv 1 \mod{9}; \\
\ls{t}{\pi}_{3}^{2} & p \equiv 4 \mod{9}; \\
\ls{t}{\pi}_{3} & p \equiv 7 \mod{9} \end{cases}$$
of \cite{prw}.  The case $m=6$ follows immediately from the fact
that the order of $\ls{a_{p}(E^{-3}_{t})}{p}_{6}$ 
equals the product of the orders of
$\ls{a_{p}(E^{-3}_{t})}{p}_{2}$ and $\ls{a_{p}(E^{-3}_{t})}{p}_{3}$.
\end{proof}

\begin{remark}
One can of course explicate the formula for $\delta_{6}^{n}(E^{-3}_{t};\C)$.
We give the formula for $n=1$ to make it clear why we do not give it in
general:
$$\delta_{6}^{1}(E^{-3}_{t};\C) = 
\rho_{\C}(\C_{36}^{1} \cap \Ctt{4}{-3}{+}) +
\rho_{\C}(\C_{36}^{19} \cap \Ct{t}{-}) +
\rho_{\C}(\C_{36}^{7,31} \cap \Ctt{6}{t}{-}) +
\rho_{\C}(\C_{36}^{13,25} \cap \Ctt{4}{-3}{+} \cap \Ctt{3}{t}{+}).$$
\end{remark}

\subsection{Abelian densities} \label{sec:ad}

For later use we now use the results of the previous section to
compute the densities
$\delta_{m'}^{1}(E^{d}_{t};\C_{m}^{1})$ for $m'$ dividing $w_{d}$ and $m$ with
$(m,w_{d})=m'$.

\begin{lemma} \label{lemma:2}
Assume that $d \neq -4$ and fix $m \geq 1$ even.  Fix
$t \in \Q^{\times}$ and let $t'$ denote the unique squarefree integer
with $t/t' \in \Q^{\times 2}$.  Then
$$\delta_{2}^{1}(E^{d}_{t};\C_{m}^{1}) = \begin{cases}
\frac{3}{4} & 4 \nmid m,\, t' \mid md,\, t' \equiv 3 \mod{4} \text{~or~}
d=-8,\, t' \equiv 2 \mod{8}; \\
\frac{1}{4} & 4 \nmid m,\, t' \mid md,\, t' \equiv 1 \mod{4} \text{~or~}
d=-8,\, t' \equiv 6 \mod{8}; \\
\frac{1}{2} & \text{otherwise}.\end{cases}$$
\end{lemma}
\begin{proof}
We assume $d \neq -8$; the case $d=-8$ is similar.
Set $\C = \C_{m}^{1} \cap 
\Ct{d}{+}$.  By Proposition~\ref{prop:m2}
we must compute the densities $\rho_{\C}(\Ctt{4}{d}{+})$ and
$\rho_{\C}(\C_{4}^{3} \cap \Ct{t}{-})$.  This is straightforward
using Lemma~\ref{lemma:dens}.  Indeed, we have
$\C = \C_{\Q(\zeta_{m},\sqrt{d})}^{\{1\}}$ and
$\Ctt{4}{d}{+} = \C_{\Q(i,\sqrt[4]{d})}^{\{1\}}$.
(In each case $\{1\}$ stands for the identity element of the corresponding
Galois group.)  Since
$$\Q(\zeta_{m},\sqrt{d}) \cap \Q(i,\sqrt[4]{d}) =
\begin{cases}
\Q(\sqrt{d}) & 4 \nmid m; \\
\Q(i,\sqrt{d}) & 4 \mid m;
\end{cases}$$
it follows from Lemma~\ref{lemma:dens} that
$$\rho_{\C}(\Ctt{4}{d}{+}) = \begin{cases}
\frac{1}{4} & 4 \nmid m; \\
\frac{1}{2} & 4 \mid m. \end{cases}$$

On the other hand, clearly $\rho_{\C}(\C_{4}^{3} \cap \Ct{t}{-})=0$
if $4$ divides $m$ or if $t'=1$.  When $4 \nmid m$ and $t' \neq 1$, 
we have
$$\C_{4}^{3} \cap \Ct{t}{-} = \C_{\Q(i,\sqrt{t})}^{\{\sigma\}}$$
with $\sigma(i)=-i$ and $\sigma(\sqrt{t})=-\sqrt{t}$.
Thus by Lemma~\ref{lemma:dens}
$$\rho_{\C}(\C_{4}^{3} \cap \Ct{t}{-}) =
\begin{cases} 0 & \sigma|_{K} \neq 1; \\
\frac{1}{[\Q(i,\sqrt{t}):K]} &
\sigma|_{K} = 1; \end{cases}$$
where
$$K = \Q(i,\sqrt{t}) \cap \Q(\zeta_{m},\sqrt{d}) =
\begin{cases} 
\Q(\sqrt{t}) & t' \mid md,\, t' \equiv 1 \mod{4}; \\
\Q(\sqrt{-t}) & t' \mid md,\, t' \equiv 3 \mod{4}; \\
\Q & \text{otherwise}. \end{cases}$$
When $K=\Q$ we thus have
$\rho_{\C}(\C_{4}^{3} \cap \Ct{t}{-})=\frac{1}{4}$.  When
$K=\Q(\sqrt{t})$ (resp.\ $K=\Q(\sqrt{-t})$) we have
$\sigma|_{K} \neq 1$ (resp.\ $\sigma|_{K} = 1$); the lemma follows
easily from this.
\end{proof}

\begin{lemma} \label{lemma:4}
Fix $m \geq 1$ even.  Let $t$ be a fourth-power free
integer and let $t'$ denote the unique squarefree integer with
$t/t'$ a square.  If $4 \nmid m$, then
$$\delta_{2}^{1}(E^{-4}_{t};\C_{m}^{1}) = \begin{cases}
1 & t' \mid m,\, t' \text{~even}; \\
\frac{1}{2} & t' \mid m,\, t' \text{~odd}; \\
\frac{3}{4} & otherwise. \end{cases}$$
If $4 \mid m$, then
$$\delta_{4}^{1}(E^{-4}_{t};\C_{m}^{1}) = \begin{cases}
\frac{3}{4} & 8 \nmid m,\, t \in \{2,-8\}; \\
\frac{1}{2} & 8 \mid m \text{~or~} t' \mid m,\, t' \text{~even},\,
t \notin \{\pm 2,\pm 8\}; \\
\frac{1}{4} & 8 \nmid m,\, t' \mid m,\, t' \text{~odd} \text{~or~}
t \in \{-2,8\}; \\
\frac{3}{8} & \text{otherwise}. \end{cases}$$
\end{lemma}
\begin{proof}
The case $4 \nmid m$ is similar to Lemma~\ref{lemma:2}; we omit the details.
Assume therefore that $4 \mid m$; we must compute the four densities
occurring in Proposition~\ref{prop:m4}.  We begin with 
$\rho_{\C_{m}^{1}}(\G_{8}^{1})$.
We have $\C^{1}_{m} = \C_{\Q(\zeta_{m})}^{\{1\}}$ and $\G_{8}^{1} =
\C_{K}^{\{1\}}$, where $K = \Q(\zeta_{16},\sqrt[4]{-2})$ is the ray class
field of $\Q(i)$ of conductor $8$.  (Note that $[K:\Q]=16$.)
As
$$\Q(\zeta_{m}) \cap K = \Q(\zeta_{(m,16)}),$$
it follows from Lemma~\ref{lemma:dens} that
$$\rho_{\C_{m}^{1}}(\G_{8}^{1}) = \frac{(m,16)}{32}.$$

We consider the two densities
$\rho_{\C_{m}^{1}}(\G_{8}^{5} \cap \Ct{t}{-})$ and
$\rho_{\C_{m}^{1}}(\G_{8}^{5 + 4i} \cap \Ct{t}{+})$ together.
These densities are both zero if $16 \mid m$, so we may assume that
$16 \nmid m$.  It follows easily from the fact that
$K$ is a non-abelian extension of $\Q$ and a quadratic extension of
$\Q(\zeta_{16})$ that
\begin{multline*}
\rho_{\C_{m}^{1}}(\G_{8}^{5} \cap \Ct{t}{-}) +
\rho_{\C_{m}^{1}}(\G_{8}^{5 + 4i} \cap \Ct{t}{+}) =
\frac{1}{2} \rho_{\C_{m}^{1}}(\C_{16}^{9} \cap \Ct{t}{-}) +
\frac{1}{2} \rho_{\C_{m}^{1}}(\C_{16}^{9} \cap \Ct{t}{+}) \\
= \frac{1}{2} \rho_{\C_{m}^{1}}(\C_{16}^{9}) = \frac{(m,8)}{32}.
\end{multline*}

It remains to compute $\rho_{\C_{m}^{1}}(\C_{8}^{5} \cap \Ctt{4}{t/2}{+})$.
This density is zero if $8 \mid m$, so we may assume that $8 \nmid m$.
We have 
$\C_{\Q(\zeta_{8},\sqrt[4]{t/2})}^{\{\sigma\}}$, where
$\sigma(\zeta_{8})=\zeta_{8}^{5}$ and $\sigma(\sqrt[4]{t/2})=\sqrt[4]{t/2}$
so long as such a $\sigma$ exists.  There is no such $\sigma$ exactly
when $t' = \pm 1 $ or $t \in \{2,-8\}$, in which case the desired density
is zero.  Otherwise, as
$$\Q(\zeta_{m}) \cap \Q(\zeta_{8},\sqrt[4]{t/2}) =
\begin{cases} \Q(i,\sqrt{t'}) & t' \mid m,\, t' \text{~odd}; \\
\Q(i,\sqrt{t'/2}) & t' \mid m,\, t' \text{~even}; \\
\Q(i) & t' \nmid m; \end{cases}$$
and $\sigma(\sqrt{t'/2})=\sqrt{t'/2}$ and $\sigma(\sqrt{2})=-\sqrt{2}$,
it is now a straightforward computation with
Lemma~\ref{lemma:dens} to determine that for $8 \nmid m$ we have
$$\rho_{\C_{m}^{1}}(\C_{8}^{5} \cap \Ctt{4}{t/2}{+}) =
\begin{cases} \frac{1}{2} & t \in \{ 2,-8\}; \\
\frac{1}{4} & t' \mid m,\, t' \text{~even}; \\
\frac{1}{8} & t' \nmid m; \\
0 & t' \mid m,\, t' \text{~odd or~} t \in \{-2,8\}.
\end{cases}
$$
The lemma follows on combining these density computations.
\end{proof}

\begin{lemma} \label{lemma:6}
Fix $m \geq 1$ divisible by $3$.  
Let $t$ be a cube-free
integer.  Then
$$\delta_{3}^{1}(E^{-3}_{t};\C_{m}^{1}) = \begin{cases} 1 &
9 \mid m \text{~or~} t=1; \\
\frac{5}{9} & 9 \nmid m \text{~and~} t \neq 1;
\end{cases}$$
\end{lemma}
\begin{proof}
This is immediate from Proposition~\ref{prop:m36} and the fact that
$$\rho_{\C_{m}^{1}}(\C_{9}^{4,7} \cap \Ctt{3}{t}{+}) =
\begin{cases} \frac{2}{3} & 9 \nmid m \text{~and~} t = 1; \\
\frac{2}{9} & 9 \nmid m \text{~and~} t \neq 1; \\
0 & 9 \mid m. 
\end{cases}$$
\end{proof}

\section{Conjectures}

\subsection{Statements}

We continue with the notation of the previous section.
Our basic conjecture is that the power residue symbol
$\ls{a_{p}(E^{d}_{t})}{p}_{m}$ is ``random'' except to the extent that
its $m/(m,w_{d})^{\th}$ power is determined by the formulae of the previous
section.  More precisely, we make the following conjecture.

\begin{conjecture} \label{conj:main}
Fix $m \geq 1$ and let $\C$ be a Chebotarev set contained in
$\C_{m}^{1} \cap \Ct{d}{+}$.  Then
$$\delta_{m}^{n}(E^{d}_{t};\C) = \frac{\pn}{m} \cdot
\frac{\delta_{m'}^{n'}(E^{d}_{t};\C)}{\varphi(n')/m'}$$
where $m'=(m,w_{d})$ and $n' = \frac{m'}{(\frac{m}{n},w_{d})}$.
\end{conjecture}

Conjecture~\ref{conj:main} is of course based on a great deal of numerical
evidence.  As the number of cases involved is somewhat overwhelming, we will
not report any of it here; we will leave it to the curious reader to
numerically verify our conjectures in any particular case.

It is worth noting that this conjecture satisfies certain compatibilities.

\begin{proposition} ~
\begin{enumerate}
\item Let $m_{1}, m_{2}$ be relatively prime integers and let
$\C$ be a Chebotarev set contained in $\Ct{d}{+}$.  
Fix divisors $n_{i}$ of each $m_{i}$ and
$t \in \Q^{\times}$.  If Conjecture~\ref{conj:main} holds for each
$\delta_{m_{i}}^{n_{i}}(E^{d}_{t};\C \cap \C_{m_{i}}^{1})$, then it holds
for $\delta_{m_{1}m_{2}}^{n_{1}n_{2}}(E^{d}_{t};\C \cap \C_{m_{1}m_{2}}^{1})$.
\item
Fix $m \geq 1$.  If Conjecture~\ref{conj:main}
holds for $\delta_{m}^{n}(E_{1};\C)$ (for all divisors $n$ of $m$ and
all Chebotarev sets $\C \subseteq \C_{m}^{1} \cap \Ct{d}{+}$), then it
holds for $\delta_{m}^{n}(E^{d}_{t};\C)$ for all $t \in \Q^{\times}$.
\end{enumerate}
\end{proposition}
\begin{proof}
The proof of (1) is straightforward; we leave it to the reader.  For
(2) we consider only the case $d \neq -3,-4$; the excluded cases are
similar but much more painful.  The case of $m$ odd is clear so we assume
that $m$ is even.  Fix $n$ dividing $m$ and a Chebotarev set
contained in $\C_{m}^{1} \cup \Ct{d}{+}$.  By (\ref{eq:tw}) we have
$$\ls{a_{p}(E^{d}_{t})}{\p}_{m} = \ls{t}{p}^{\frac{p-1}{m}}
\ls{a_{p}(E^{d}_{1})}{\p}_{m}$$
where $\p$ is some prime of $\Q(\zeta_{m})$ lying over $p$.
It follows that
$\ls{a_{p}(E^{d}_{t})}{\pi}_{m}$ and $\ls{a_{p}(E^{d}_{1})}{\pi}_{m}$ have
the same order unless $\frac{p-1}{m}$ is odd,
$\ls{t}{p}=-1$ and one has odd order, in which case the other has twice
that order.  That is,
\begin{equation} \label{eq:pmn}
\P_{m}^{n}(E^{d}_{t};\C) = \P_{m}^{n}(E^{d}_{1};\C \cap \C^{1}_{2m}) \cup
\P_{m}^{n}(E^{d}_{1};\C \cap \C^{m+1}_{2m} \cap \Ct{t}{+}) \cup
\P_{m}^{\tilde{n}}(E^{d}_{1};\C \cap \C^{m+1}_{2m} \cap \Ct{t}{-})
\end{equation}
where
$$\tilde{n} = \begin{cases} 2n & n \text{~odd}; \\
n/2 & n \equiv 2 \mod{4}; \\
n & 4 \mid n. \end{cases}$$

When $m$ is divisible by $4$, we either have that $n$ is also divisible
by $4$ or else $n'=1$.  It follows easily that 
applying Conjecture~\ref{conj:main}
and Proposition~\ref{prop:m2} to (\ref{eq:pmn}) thus yields
\begin{align*}
\delta_{m}^{n}(E^{d}_{t};\C) &=
\frac{2\varphi(n)}{m} \cdot \left(
\rho_{\C \cap \C_{2m}^{1}}(\Ctt{4}{d}{+}) +
\rho_{\C \cap \C_{2m}^{m+1} \cap \Ct{t}{+}}(\Ctt{4}{d}{+}) +
\rho_{\C \cap \C_{2m}^{m+1} \cap \Ct{t}{-}}(\Ctt{4}{d}{+}) \right) \\
&= 2\frac{\varphi(n)}{m} \cdot \rho_{\C}(\Ctt{4}{d}{+}) \\
&= \frac{\varphi(n)}{m} \cdot \frac{\delta_{2}^{n'}(E^{d}_{t};\C)}
{\varphi(n')/m'}
\end{align*}
as desired.

When $m \equiv 2 \mod{4}$ we have
$\C \cap \C_{2m}^{1} = \C \cap \C_{4}^{1}$ and
$\C \cap \C_{2m}^{m+1} = \C \cap \C_{4}^{3}$.  
If $n$ is even, so that $n'=2$, $\tilde{n}=\frac{n}{2}$ and 
$\tilde{n}' = 1$, from (\ref{eq:pmn}) we obtain
\begin{align*}
\delta_{m}^{n}(E^{d}_{t};\C) &=
\frac{2\varphi(n)}{m} \cdot \left(
\rho_{\C \cap \C_{4}^{1}}(\Ctt{4}{d}{-}) +
\rho_{\C \cap \C_{4}^{3} \cap \Ct{t}{+}}(\C_{4}^{3}) +
\rho_{\C \cap \C_{4}^{3} \cap \Ct{t}{-}}(\emptyset) \right) \\
&= \frac{2\varphi(n)}{m} \cdot \left(
\rho_{\C \cap \C_{4}^{1}}(\Ctt{4}{d}{-}) +
\rho_{\C}(\C_{4}^{3} \cap \Ct{t}{+})  \right) \\
&= \frac{\varphi(n)}{m} \cdot \frac{\delta_{2}^{n'}(E^{d}_{t};\C)}
{\varphi(n')/m'}.
\end{align*}
Finally, if $n$ is odd we have $n'=1$, $\tilde{n} = 2n$ and $\tilde{n}=2$, so
that (\ref{eq:pmn}) yields
\begin{align*}
\delta_{m}^{n}(E^{d}_{t};\C) &=
\frac{2\varphi(n)}{m} \cdot \left(
\rho_{\C \cap \C_{4}^{1}}(\Ctt{4}{d}{+}) +
\rho_{\C \cap \C_{4}^{3} \cap \Ct{t}{+}}(\emptyset) +
\rho_{\C \cap \C_{4}^{3} \cap \Ct{t}{-}}(\C_{4}^{3}) \right) \\
&= \frac{2\varphi(n)}{m} \cdot \left(
\rho_{\C \cap \C_{4}^{1}}(\Ctt{4}{d}{+}) +
\rho_{\C}(\C_{4}^{3} \cap \Ct{t}{-})  \right) \\
&= \frac{\varphi(n)}{m} \cdot \frac{\delta_{2}^{n'}(E^{d}_{t};\C)}
{\varphi(n')/m'}.
\end{align*}
\end{proof}

\subsection{Abelian densities}

Combining Conjecture~\ref{conj:main} with the calculations of 
Section~\ref{sec:ad} we obtain the following explicit conjectures for
the densities $\delta_{m}^{1}(E^{d}_{t};\C_{m}^{1})$.

\begin{proposition} \label{prop:abc}
Fix $m \geq 1$ and fix a non-zero integer $t$; if $d=-4$ (resp.\ $d=-3$)
assume also that $t$ is fourth-power free (resp.\ sixth-power-free).
Let $t'$ denote the unique
squarefree integer with $t/t' \in \Q^{\times 2}$. Assume that
Conjecture~\ref{conj:main} holds for $E^{d}_{t}$ and $m$.
\begin{enumerate}
\item If $(m,w_{d})=1$, then
$$\delta_{m}^{1}(E^{d}_{t};\C_{m}^{1} \cap \Ct{d}{+}) = \frac{1}{m}.$$
\item If $(m,w_{d})=2$ and $d \neq -4$,  then
$$\delta_{m}^{1}(E^{d}_{t};\C_{m}^{1} \cap \Ct{d}{+}) =
\begin{cases} \frac{3}{2m} & 4 \nmid m,\, t' \mid md,\, t' \equiv 3 \mod{4}
\\ & \text{~or~} d=-8,\, t' \equiv 2 \mod{8}; \\
\frac{1}{2m} & 4 \nmid m,\, t' \mid md,\, t' \equiv 1 \mod{4} \\
& \text{~or~} d=-8,\, t' \equiv 6 \mod{8}; \\
\frac{1}{m} & \text{otherwise.} \end{cases}$$
\item If $(m,w_{d})=2$ and $d=-4$, then
$$\delta_{m}^{1}(E^{d}_{t};\C^{1}_{m} \cap \C^{1}_{4}) =
\begin{cases}
\frac{2}{m} & t' \mid m,\, t' \text{~even}; \\
\frac{1}{m} & t' \mid m,\, t' \text{~odd}; \\
\frac{3}{2m} & \text{otherwise.} \end{cases}$$
\item If $(m,w_{d})=3$ (so that $d=-3$), then
$$\delta_{m}^{1}(E^{d}_{t};\C_{m}^{1}) =
\begin{cases} \frac{3}{m} & 9 \mid m \text{~or~} t \in \Q^{\times 3}; \\
\frac{5}{3m} & \text{otherwise.}
\end{cases}$$
\item If $(m,w_{d})=4$ (so that $d=-4)$, then
$$\delta_{m}^{1}(E^{d}_{t};\C_{m}^{1}) =
\begin{cases} 
\frac{3}{m} & 8 \nmid m,\, t \in \{2,-8\}; \\
\frac{2}{m} & 8 \mid m \text{~or~} t' \mid m,\, t' \text{~even~},\, 
t \notin \{ \pm 2,\pm 8\}; \\
\frac{1}{m} & 8 \nmid m,\, t' \mid m,\, t' \text{~odd or~} t \in \{-2,8\}; \\
\frac{3}{2m} & \text{otherwise.} \end{cases}$$
\item If $(m,w_{d})=6$ (so that $d=-3)$, then
$$\delta_{m}^{1}(E^{d}_{t};\C_{m}^{1}) = \begin{cases} 
\frac{9}{2m} & 4 \nmid m,\, t' \mid md,\, t' \equiv 3 \mod{4},\, t \in \Q^{\times 3};\\
\frac{3}{2m} & 4 \nmid m,\, t' \mid md,\, t' \equiv 1 \mod{4},\, t \in \Q^{\times 3}; \\
\frac{3}{m} & 4 \mid m \text{~or~} t' \nmid md \text{~or~} t' \text{~even
and~} t \in \Q^{\times 3}; \\
\frac{5}{2m} & 4 \nmid m,\, t' \mid md,\, t' \equiv 3 \mod{4},\, t \notin \Q^{\times 3};\\
\frac{5}{6m} & 4 \nmid m,\, t' \mid md,\, t' \equiv 1 \mod{4},\, t \notin \Q^{\times 3}; \\
\frac{5}{3m} & 4 \mid m \text{~or~} t' \nmid md \text{~or~} t' \text{~even
and~} t \notin \Q^{\times 3};
\end{cases}$$
\end{enumerate}
\end{proposition}
\begin{proof}
This is all immediate from the results of Section~\ref{sec:ad}.
\end{proof}

\section{Biquadratic residues} \label{sec:gi}

In this section we give the proof of Proposition~\ref{prop:m4}.  Recall
that in this setting
$E^{-4}_{t}$ denotes the elliptic curve $y^{2} = x^{3} - tx$ with
complex multiplication by $\Z[i]$.  By \cite{Silverman2} we have
that $a_{p}(E^{-4}_{t})=0$ for $p \equiv 3 \mod{4}$ or not relatively
prime to $t$; otherwise
\begin{equation} \label{eq:zi}
a_{p}(E^{-4}_{t}) = \ls{t}{\pib}_{4}\pi + \ls{t}{\pi}_{4}\pib
\end{equation}
where $p=\pi\pib$ with $\pi,\pib$ primary irreducibles in $\Z[i]$.

We begin with a computation with biquadratic reciprocity.

\begin{lemma} \label{lemma:div}
Let $\pi = a+bi$ be a primary irreducible of prime norm $p$
in $\Z[i]$ and let $\ell$
be an odd rational prime divisor of $a$.  Then
$$\ls{\ell}{\pi}_{4} = (-1)^{\frac{\ell-1}{2}\cdot\frac{p-1}{4}}\ls{2}{\ell}.$$
\end{lemma}
\begin{proof}
If $\ell \equiv 1 \mod{4}$, then $\ell$ factors as $\ell = \lambda\bar{\lambda}$
into primary irreducibles.  Thus
$$\ls{\ell}{\pi}_{4} = \ls{\lambda}{\pi}_{4}\ls{\bar{\lambda}}{\pi}_{4} =
\ls{\pi}{\lambda}_{4}(-1)^{\frac{\ell-1}{4}\cdot \frac{p-1}{4}} \cdot
\ls{\pi}{\bar{\lambda}}_{4}(-1)^{\frac{\ell-1}{4}\cdot \frac{p-1}{4}} =
\ls{\pi}{\lambda}_{4}\ls{\pi}{\bar{\lambda}}_{4}$$
by biquadratic reciprocity.  Thus
$$\ls{\ell}{\pi}_{4} = \ls{a+bi}{\lambda}_{4}\ls{a+bi}{\bar{\lambda}}_{4} =
\ls{bi}{\lambda}_{4}\ls{bi}{\bar{\lambda}}_{4}$$
since $\ell$ divides $a$.  As
$$\overline{\ls{bi}{\lambda}}_{4} = \ls{-bi}{\bar{\lambda}}_{4}$$
we conclude that
$$\ls{\ell}{\pi}_{4} = \ls{-1}{\lambda}_{4} = (-1)^{\frac{\ell-1}{4}}.$$
Since
$$\frac{\ell-1}{4} \equiv \frac{\ell^{2}-1}{8} \pmod{2}$$
for $\ell \equiv 1 \pmod{4}$, the lemma follows in this case.

If $\ell \equiv 3 \mod{4}$, then $-\ell$ is primary,
so that by biquadratic reciprocity we have
$$\ls{\ell}{\pi}_{4} = \ls{-1}{\pi}_{4}\ls{-\ell}{\pi}_{4} =
(-1)^{\frac{p-1}{4}}
 \ls{\pi}{\ell}_{4}(-1)^{\frac{\ell^{2}-1}{4} \cdot
\frac{p-1}{4}} = (-1)^{\frac{p-1}{4}}\ls{\pi}{\ell}_{4}$$
since $\ell^{2}-1$ is divisible by $8$.  We now have
$$\ls{\ell}{\pi}_{4} = (-1)^{\frac{p-1}{4}}
\ls{\pi}{\ell}_{4} = (-1)^{\frac{p-1}{4}}\ls{bi}{\ell} = (-1)^{\frac{p-1}{4}}
\ls{b}{\ell}_{4}
i^{\frac{\ell^{2}-1}{4}}.$$
Since
$$b^{\frac{\ell^{2}-1}{4}} \equiv 1 \pmod{\ell}$$
by Fermat's little theorem, the lemma follows in this case from the fact
that
$$i^{\frac{\ell^{2}-1}{4}} = (-1)^{\frac{\ell^{2}-1}{8}} = \ls{2}{\ell}.$$
\end{proof}

It will also be useful to recall the computation of the
biquadratic character of $2$.

\begin{lemma} \label{lemma:bq2}
Let $a+bi$ be a primary irreducible of norm $p$.  Then
$$\ls{2}{a+bi}_{4} = (-i)^{b/2}.$$
\end{lemma}
\begin{proof}
Factoring $2$ as
$$2 = -i(1+i)^{2},$$
we find that
$$\ls{2}{a+bi}_{4} = \ls{-i}{\pi}_{4} \cdot \ls{1+i}{\pi}_{4}^{2} =
(-i)^{\frac{p-1}{4}} \cdot i^{\frac{a-b-b^{2}-1}{2}}$$
by \cite[Theorem 6.9]{lemmermeyer}.  From here a lengthy but
elementary calculation yields the asserted formula.
\end{proof}

We are now in a position to give a formula for
the biquadratic residue symbol $\ls{a_{p}(E^{-4}_{t})}{\pi}_{4}$ with
$\pi$ a primary irreducible divisor of $p$.

\begin{proposition} \label{prop:api}
Let $p \equiv 1 \mod{4}$ be a prime relatively prime to $t$ and let
$a+bi$ be a primary irreducible of norm $p$.  Then
$$\ls{a_{p}(E^{-4}_{t})}{a+bi}_{4} = i^{\frac{a-1}{2}} \cdot
(-1)^{\frac{b^{2}+2b}{8}} \cdot
\ls{t}{a+bi}_{4}^{\frac{p-1}{4}}.$$
\end{proposition}
\begin{proof}
Set $\pi = a+bi$.
By (\ref{eq:zi}) we have
$$a_{p}(E^{-4}_{t}) \equiv \ls{t}{\pi}_{4}\pib \pmod{\pi}$$
so that
$$\ls{a_{p}(E^{-4}_{t})}{\pi}_{4} =
\ls{\ls{t}{\pi}_{4}\pib}{\pi}_{4} = \ls{t}{\pi}_{4}^{\frac{p-1}{4}}
\ls{\pib}{\pi}_{4}$$
It thus suffices to compute $\ls{\pib}{\pi}_{4}$.

Let
$$a = p_{1}^{e_{1}} \cdots p_{r}^{e_{r}}q_{1}^{f_{1}} \cdots
q_{s}^{f_{s}}$$
be the prime of factorization of $a$, where $p_{i} \equiv 1 \mod{4}$
and $q_{i} \equiv 3 \mod{4}$.  (Here $a$ is odd since $\pi$ is primary.)
Since
$$\pib \equiv \pi + \pib \equiv 2a \pmod{\pi},$$
applying Lemma~\ref{lemma:div} we find that
\begin{align*}
\ls{\pib}{\pi}_{4} &= \ls{2}{\pi}_{4}\ls{p_{1}}{\pi}_{4}^{e_{1}}\cdots
\ls{p_{r}}{\pi}_{4}^{e_{r}}\ls{q_{1}}{\pi}_{4}^{f_{1}}\cdots
\ls{q_{s}}{\pi}_{4}^{f_{s}} \\
&= \ls{2}{\pi}_{4}\ls{2}{p_{1}}^{e_{1}}\cdots \ls{2}{p_{r}}^{e_{r}}
\ls{2}{q_{1}}^{f_{1}}\cdots \ls{2}{q_{s}}^{f_{s}}
(-1)^{\frac{p-1}{4} \cdot (f_{1}+\cdots+f_{s})}.
\end{align*}
By the multiplicativity of the Jacobi symbol and the fact that
$$f_{1} + \cdots + f_{s} \equiv \frac{a-1}{2} \pmod{2},$$
we conclude that
$$\ls{\pib}{\pi}_{4} = \ls{2}{\pi}_{4}
\ls{2}{a}(-1)^{\frac{p-1}{4} \cdot \frac{a-1}{2}}.$$
Combining this with Lemma~\ref{lemma:bq2} and simplifying yields the
theorem.
\end{proof}

Proposition~\ref{prop:m4} follows immediately.

\begin{corollary}
Let $\C$ be a Chebotarev set.  For any $t \in \Q^{\times}$ we have
\begin{align*}
\delta_{4}^{1}(E^{-4}_{t};\C) &= \rho_{\C}(\G_{8}^{1}) +
\rho_{\C}(\G_{8}^{5} \cap \Ct{t}{-}) +
\rho_{\C}(\G_{8}^{5+4i} \cap \Ct{t}{+}) +
\rho_{\C}(\C_{8}^{5} \cap \Ctt{4}{t/2}{+}); \\
\delta_{4}^{2}(E^{-4}_{t};\C) &= \rho_{\C}(\G_{8}^{1+4i}) +
\rho_{\C}(\G_{8}^{5} \cap \Ct{t}{+}) +
\rho_{\C}(\G_{8}^{5+4i} \cap \Ct{t}{-}) +
\rho_{\C}(\C_{8}^{5} \cap \Ctt{4}{t/2}{-}); \\
\delta_{4}^{4}(E^{-4}_{t};\C) &= \rho_{\C}(\C_{8}^{5} \cap \Ct{t}{+}).
\end{align*}
\end{corollary}
\begin{proof}
Since the nine sets of primes listed above partition the set of all
primes (relatively prime to $t$), it suffices to show that
$\ls{a_{p}(E^{-4}_{t})}{\pi}_{4}$ has the asserted order for each set.
This is straightforward from Proposition~\ref{prop:api}.
Indeed, fix $p \in \C_{4}^{1}$ and a primary divisor $\pi$ of $p$.
Proposition~\ref{prop:api} (together with Lemma~\ref{lemma:bq2} when
$p \equiv 5 \mod{8}$) then yields
$$\ls{a_{p}(E^{-4}_{t})}{\pi}_{4} = \begin{cases}
1 & p \in \G_{8}^{1}; \\
-1 & p \in \G_{8}^{1+4i}; \\
-\ls{t}{p} & p \in \G_{8}^{5}; \\
\ls{t}{p} & p \in \G_{8}^{5+4i}; \\
\ls{t/2}{\pi}_{4} & p \in \C_{8}^{5}.
\end{cases}$$
The corollary follows easily.
\end{proof}

\end{document}